\documentclass[10pt]{amsart}
\usepackage{latexsym, amsmath, amssymb}

\renewcommand{\epsilon}{\varepsilon}
\newcommand{\newsection}[1]
{\subsection{#1}\setcounter{theorem}{0} \setcounter{equation}{0}
\par\noindent}

\newtheorem{theorem}{Theorem}

\newtheorem{lemma}[theorem]{Lemma}
\newtheorem{corr}[theorem]{Corollary}

\newtheorem{proposition}[theorem]{Proposition}

\newcommand{\cd}{\, \cdot\, }

\newcommand{\R}{{\mathbb R}}
\newcommand{\bdy}{{\partial\Omega}}
\newcommand{\T}{{T_\varepsilon}}
\newcommand{\la}{\langle}
\newcommand{\ra}{\rangle}
\newcommand{\tGamma}{\tilde{\Gamma}}

\begin{document}

\title[Quasilinear wave equations in waveguides]
{Almost global existence for quasilinear wave equations in waveguides with Neumann
boundary conditions}

\thanks{The authors were supported in part by the NSF}
\thanks{A portion of this work was completed while the authors were visiting the Mathematical
Sciences Research Institute (MSRI).  The authors gratefully acknowledge the 
hospitality and support of MSRI}

\author{Jason Metcalfe}
\address{Department of Mathematics, University of California, Berkeley, CA  94720-3840}
\email{metcalfe@math.berkeley.edu}
\author{Ann Stewart}
\address{Department of Mathematics,  Johns Hopkins University,
Baltimore, MD 21218}
\email{stewart@jhu.edu}

%\maketitle

\begin{abstract}

In this paper, we prove almost global existence of solutions to certain quasilinear wave
equations with quadratic nonlinearities in infinite homogeneous waveguides with Neumann boundary
conditions.  
We use a Galerkin method to expand the Laplacian of the compact
base in terms of its eigenfunctions.  For those terms corresponding to zero modes, we obtain decay
using analogs of estimates of Klainerman and Sideris.  For the nonzero modes, estimates
for Klein-Gordon equations, which provide better decay, are available.

\end{abstract}
\maketitle

\newsection{Introduction}

The purpose of this paper is to give a simple proof of almost global existence for quasilinear
Neumann wave equations on infinite homogeneous waveguides.  This work expands on results
of Metcalfe, Sogge, and Stewart \cite{MSS} and of Lesky and Racke \cite{lesky}.  Here, as in these works,
the key step will be the use an eigenfunction expansion in the compact base and to use
estimates for Klein-Gordon equations for those terms corresponding to nonzero modes.  For the zero modes,
we will use decay estimates that are analogous to those of Klainerman and Sideris \cite{KSid}.  In order
to prove long time existence, we couple these decay estimates with energy estimates.  

Let us describe our initial-boundary value problem more precisely.  We will be studying nonlinear 
wave equations in infinite homogeneous waveguides, $\R^3\times \Omega$ where $\Omega\subset \R^d$
denotes a nonempty, bounded domain with smooth boundary.  We examine equations of the form
\begin{equation}
\label{main.equation}
\begin{cases}
\Box u = Q(\partial_{t,x}u, \partial_{t,x}^2 u),\quad (t,x,y)\in \R_+\times\R^3\times\Omega,\\
\partial_\nu u(t,x,\cd)|_\bdy = 0,\\
u(0,x,y)=f(x,y),\quad \partial_t u(0,x,y)=g(x,y),
\end{cases}
\end{equation}
where $\Box=\partial_t^2 - (\Delta +\Delta_\Omega)$ is the d'Alembertian on the waveguide.  Here,
$\partial_\nu$ denotes the normal derivative on $\bdy$.  Moreover,
$$\Delta=\Delta_{\R^3}=\sum_{j=1}^3 \partial^2/\partial x_j^2$$
is the Laplacian on $\R^3$ and 
$$\Delta_\Omega = \sum_{j=1}^d \partial^2 / \partial y_j^2$$
is the Neumann Laplacian.  

The nonlinearity $Q$ is quadratic in its arguments and is affine linear in $\partial_{t,x}^2 u$ (i.e. is
quasilinear).  We may expand $Q$ as follows
\begin{equation}
\label{Q}
Q(\partial_{t,x} u, \partial_{t,x}^2u)=\sum_{0\le j,k,l\le 3} A^{jk}_l \partial_l u \partial_j\partial_k u
+ R(\partial_{t,x}u, \partial_{t,x}u)
\end{equation}
where the $A^{jk}_l$ are real constants and $R$ is a constant coefficient, quadratic form.  Here and throughout,
we set $x_0=t$ and $\partial_0=\partial_t$ when convenient.

In order to solve \eqref{main.equation}, the data must be assumed to satisfy the relevant compatibility
condition.  Let $J_k u = \{\partial_{x,y}^\alpha u\, :\, 0\le |\alpha|\le k\}$ denote the
collection of all spatial derivatives of $u$ of order up to $k$ (using local coordinates
in a small tubular neighborhood of $\bdy$).  If $u$ is a formal $H^N$ solution for some fixed $N$,
then we can write $\partial_t^k \partial_\nu u(0,\cd)=\Psi_k(J_{k+1}f,J_kg)$ where the 
$\Psi_k$ are called compatibility functions and depend on $Q$, $J_{k+1}f$ and $J_k g$.
The compatibility conditions for $(f,g)\in H^{N}\times H^{N-1}$ simply require that the $\Psi_k$
vanish on $\R_+\times\R^3\times \bdy$ when $0\le k\le N-2$.  Moreover, we say that $(f,g)\in C^\infty$
satisfy the compatibility conditions to infinite order if this condition holds for all $N$.

We assume that the initial data have compact support and are small in norm.  That is, we assume that
there is a fixed constant $B>0$ so that
\begin{equation}
\label{data.support}
f(x,y)=g(x,y)=0,\quad |x|>B.
\end{equation}
Moreover, we assume that 
\begin{equation}
\label{data.smallness}
\|f\|_{H^N(\R^3\times\Omega)} + \|g\|_{H^{N-1}(\R^3\times\Omega)} \le \varepsilon,
\end{equation}
where 
$$\|f\|_{H^N(\R^3\times\Omega)} = \sum_{|\alpha|\le N} \|\partial_{x,y}^\alpha f\|_{L^2(\R^3\times\Omega)}.$$

Under these assumptions, we can prove the following almost global existence result.
\begin{theorem}
\label{theorem1.1}
Assume that the Cauchy data $(f,g)\in C^\infty(\R^3\times\Omega)$ satisfy \eqref{data.support} and
\eqref{data.smallness} as well as the compatibility conditions to infinite order.  Then there are constants
$N$, $\kappa$ and $\varepsilon_0>0$ so that if $\varepsilon<\varepsilon_0$ and $N$ is sufficiently large 
in \eqref{data.smallness}, then
\eqref{main.equation} has a unique solution $u\in C^\infty([0,\T)\times \R^3\times\Omega)$ where
\begin{equation}\label{lifespan}
\T=\exp(\kappa/\varepsilon).
\end{equation}
\end{theorem}

Notice that the lifespan \eqref{lifespan} is sharp.  Indeed, if one takes $f(x,y)=f(x)$, $g(x,y)=g(x)$
independent of $y$, then solutions of \eqref{main.equation} are equivalent to solutions to
$\Box u=Q(\partial_{t,x}u, \partial^2_{t,x} u)$ in Minkowski space.  By the classical counterexamples
of John (see, e.g., \cite{john1, john2}) of the form $\Box u = (\partial_t u)^2$, it is seen
that \eqref{lifespan} cannot be improved.  See, also, Sideris \cite{Sid}.

It is for technical reasons that we are only able to handle nonlinearities that do not depend on 
the derivatives $\partial_y$.  However, when the compact base is one dimensional, we are able to obtain
an optimal result.  In order to use energy methods, we must now assume a nonlinear compatibility
condition as in \cite{MSS}.  Here, we will be studying the initial-boundary value problem
\begin{equation}
\label{main.equation.2}\begin{cases}
\Box u = \tilde{Q}(\partial u, \partial^2 u),\quad (t,x,y)\in \R_+\times \R^3\times [a,b],\\
\partial_y u(t,x,a)=0,\quad \partial_y u(t,x,b)=0\\
u(0,x,y)=f(x,y),\quad \partial_t u(0,x,y)=g(x,y),
\end{cases}
\end{equation}
where $a<b$ are fixed constants and $\partial=\partial_{t,x,y}$ denotes the full space-time gradient.
Expanding $\tilde{Q}$ as above, we have that
\begin{equation}\label{Q.2}
\tilde{Q}(\partial u,\partial^2 u)=\sum_{0\le j,k,l\le 4} B^{jk}_l \partial_l u \partial_j \partial_k u
+ \tilde{R}(\partial u, \partial u),
\end{equation}
where $\tilde{R}$ is again a constant coefficient quadratic form.  Here, for convenience,
we set $x_4 = y$ and $\partial_4 = \partial_y$.  In order to use energy methods with Neumann boundary conditions,
we must assume the following nonlinear Neumann compatibility condition
\begin{equation}
\label{Neumann.compatibility}
\sum_{0\le j,k,l\le 4} B^{jk}_l \xi_l \eta_j \theta_k = 0,\quad \text{if } (\theta,\xi,\eta)\in X,
\end{equation}
where
\begin{equation}
\label{X}
X=\{(\theta,\xi,\eta)\, : \, \theta=(0,0,0,\nu(y)),\, \xi\cdot\theta=0,\,\eta\cdot\theta=0,\, y\in \bdy\}.
\end{equation}
I.e., when $\theta$ is normal to $R^{1+3}\times\bdy$ and $\xi$ and $\eta$ are orthogonal to $\theta$,
\eqref{Neumann.compatibility} must hold.  As this condition automatically holds when the quasilinear
terms only involve $\partial_j \partial_k u$, $0\le j,k\le 3$, such an assumption was unnecessary 
when studying \eqref{main.equation}.  This is the natural condition that is required in order
to obtain energy estimates for quasilinear equation with Neumann boundary conditions.

Under these assumptions, we can then prove
\begin{theorem}
\label{theorem1.2}
Assume that the data $(f,g)\in \R^3\times [a,b]$ satisfy \eqref{data.support} and \eqref{data.smallness}
as well as the compatibility conditions to infinite order.  Moreover, assume that $\tilde{Q}$ satisfies
\eqref{Neumann.compatibility}.
Then, there are constants $N$, $\kappa$, and
$\varepsilon_0$ so that if $\varepsilon<\varepsilon_0$ and $N$ is sufficiently large in \eqref{data.smallness},
then \eqref{main.equation.2} has a unique solution $u\in C^\infty([0,\T)\times \R^3\times [a,b])$.
\end{theorem}

Results similar to Theorem \ref{theorem1.1} and Theorem \ref{theorem1.2} were first obtained by
Lesky and Racke \cite{lesky} in $\R_+\times\R^n\times\Omega$ when $n\ge 5$.  These results were improved
in Metcalfe, Sogge, and Stewart \cite{MSS}.  There, for Dirichlet boundary conditions, global existence
was shown in $\R_+\times\R^n\times \Omega$ where $\Omega\subset \R^d$ for any $n\ge 3$.  By using an eigenvalue
expansion, one can reference decay estimates for the Klein-Gordon equation which provide $O(t^{-n/2})$ decay
as opposed to the standard $O(t^{-(n-1)/2})$ decay that is available for the wave equation.

When Neumann boundary conditions are assumed, the proofs are more delicate.  First of all, in order to use
energy methods a condition of the form \eqref{Neumann.compatibility} must be assumed.  This is then
sufficient to prove global existence for Klein-Gordon equations.  For the wave equation, 
as there
is a zero mode, the estimates that are available do not provide the same decay.  In \cite{MSS}, the techniques
of \cite{KSS2} were applied to obtain existence for semilinear wave equations whose quadratic nonlinearities
only depend on $\partial_x u$.  

Theorem \ref{theorem1.1} and Theorem \ref{theorem1.2} are improvements
to this last result of \cite{MSS}.  Here, we instead use the eigenfunction expansion.  For those
terms corresponding to zero modes, we use estimates analogous to those of \cite{KSid} for the wave equation.
For the nonzero modes, we have the better estimates for Klein-Gordon equations.  As the estimates of
\cite{KSid} only provide decay for the gradient $\partial_{t,x} u$, we are lead to the restriction
on the nonlinearity in Theorem \ref{theorem1.1}.  Theorem \ref{theorem1.2} follows from the simple observation
that $\partial_y u=\partial_\nu u$ when the compact base is one dimensional.  Thus, in this case, $\partial_y u$
solves a wave equation with Dirichlet boundary conditions.  Hence, better decay estimates are 
available for these terms in this special case.

This paper is organized as follows.  In the next section, we gather some weighted Sobolev estimates and
estimates for the wave equation that are analogs of those in \cite{KSid}.  In the following
section, we show that our necessary decay estimate for the wave equation follows from those estimates
in the previous section.  We also provide the decay estimate that we shall use for solutions to Klein-Gordon
equations.  This is essentially from \cite{H}.  In the fourth section, we extend the previous estimates
(which were on $\R_+\times\R^3$) to estimates on the waveguide.  To do so, we use a Galerkin expansion
in the $y$ variable and use the inequalities from the previous section.  In the fifth section, we examine
the energy estimates which we shall require.  These are rather standard, except that we need to make a natural
nonlinear compatibility assumption on the quasilinear terms.  Such an assumption is necessitated by
the Neumann boundary conditions.  Finally, in the last section, we prove Theorem \ref{theorem1.1} and
Theorem \ref{theorem1.2}.

It is a pleasure to thank C. Sogge for his collaboration \cite{MSS} that preceded this paper 
and for helpful suggestions concerning the current study.

%%%%%%%%%%%%%%%%%%%%%%%%%%%%%%%%%%%%%%%%%%%%%%%%%%%%%%%%%%%%%%%%%%%%%%%%%%%%%%%%%%%%%%%%%%%%%%%%%%%%%%
%%%%%%%%%%%%%%%%%%%%%%%%%%%%%%%%%%%%%%%%%%%%%%%%%%%%%%%%%%%%%%%%%%%%%%%%%%%%%%%%%%%%%%%%%%%%%%%%%%%%%%
\newsection{Sobolev estimates and decay estimates for wave equations in $\R_+\times\R^3$}

In this section, we present the estimates that provide the necessary decay for wave equations in
$\R_+\times\R^3$.  These will be used to estimate the terms in the expansion of $u$ with zero mode.
Here, we let
\begin{equation}
\label{gamma}
\{\Gamma\}=\{\partial_t, \partial_x, \Omega_{jk}\,:\, 0\le j<k\le 3\}
\end{equation}
where
$$\Omega_{ij}=x_i\partial_j - x_j\partial_i, \quad 1\le i<j\le 3$$
and
$$\Omega_{0k}=x_k \partial_t + t\partial_k,\quad 1\le k\le 3.$$
We will also denote
$$Z=\{\partial_t, \partial_x, \Omega_{jk}\,:\, 1\le j<k\le 3\}.$$

The first lemma that we shall require is the following weighted Sobolev inequality.
\begin{lemma}
\label{lemma.weighted.Sobolev}
Let $w\in C^\infty_0(\R_+\times\R^3)$.  Then,
\begin{equation}\label{sobolev.1}
\la r\ra^{1/2} |w(t,x)|\lesssim \sum_{|\alpha|\le 1} \|Z^\alpha \partial w(t,\cd)\|_2,
\end{equation}
and
\begin{equation}
\label{sobolev.2}
\la r\ra |\partial \Gamma^\alpha w(t,x)|\lesssim \sum_{|\beta|\le |\alpha|+2} \|\Gamma^\beta \partial w(t,\cd)\|_2.
\end{equation}
\end{lemma}

Here and throughout, we use the notation $A\lesssim B$ to indicate that $A\le CB$ for some positive, unspecified
constant $C$.  We are also using $\la x\ra = \la r\ra=\sqrt{1+|x|^2}$.
The first of these estimates essentially appears in Sideris \cite{Si}.  
See also Hidano \cite{Hid}.  The second estimate
is from Klainerman and Sideris \cite{KSid}.  As a corollary of \eqref{sobolev.1}, we obtain
\begin{corr}
\label{corollary.hidano}
Let $w\in C^\infty_0(\R_+\times\R^3)$.  Then,
\begin{equation}
\label{hidano.1}
\la r\ra^{1/2} \la t-r\ra |\partial \Gamma^\alpha w(t,x)|\lesssim \sum_{|\beta|\le |\alpha|+1}
\|\Gamma^\beta w'(t,\cd)\|_2 + \sum_{|\beta|\le |\alpha|+1} \|\la t-r\ra \partial^2 \Gamma^\beta w(t,\cd)\|_2.
\end{equation}
\end{corr}

This estimate also appeared in \cite{Hid} and follows by simply applying \eqref{sobolev.1} to 
$\langle t-r\rangle \partial \Gamma^\alpha w(t,x)$.

In order to estimate the last term in \eqref{hidano.1}, we will require estimates that are analogs of those
due to Klainerman and Sideris \cite{KSid}.  In \cite{KSid}, the necessary estimates were restricted to only make
use of the vector fields $Z$ and the scaling vector field $L=t\partial_t + r\partial_r$.  In the current
setting, it is more convenient to use estimates that only involve $\Gamma$.  Thus, we prove the following
variant of the estimates of \cite{KSid}.
\begin{lemma}
\label{lemma.KS.estimates}
Let $w\in C^2(\R_+\times\R^3)$.  Then,
\begin{align}
\la t-r\ra |\Delta w(t,x)|&\lesssim \sum_{|\alpha|\le 1} |\partial \Gamma^\alpha w(t,x)|
+(t+r)|\Box w(t,x)|,\label{KS.delta}\\
\la t-r\ra |\partial_t^2 w(t,x)|&\lesssim \sum_{|\alpha|\le 1} |\partial \Gamma^\alpha w(t,x)|
+ t|\Box w(t,x)|,\label{KS.dtdt}\\
\la t-r\ra |\nabla_x \partial_t w(t,x)|&\lesssim \sum_{|\alpha|\le 1} |\partial \Gamma^\alpha w(t,x)|
+t|\Box w(t,x)|.\label{KS.dtdx}
\end{align}
\end{lemma}
Here and throughout the current section and the next, we are abusing notation and 
using $\Box=\partial_t^2-\Delta_{\R^3}$ to denote the d'Alembertian
on Minkowski space.

The proof is based in part on the following lemma, which also appeared in \cite{KSid}.
\begin{lemma}
\label{lemma.Delta.drdr}
Let $w\in C^2(\R^3)$.  Then,
\begin{equation}
\label{Delta.drdr}
|\Delta w(x)-\partial_r^2 w(x)|\lesssim \frac{1}{r}\sum_{|\alpha|\le 1} |\nabla_x Z^\alpha w(x)|.
\end{equation}
\end{lemma}

This estimate follows immediately from the expansion of the Laplacian into its radial and
angular parts
$$\Delta=\partial_r^2 + \frac{2}{r}\partial_r + \frac{1}{r^2}\Omega\cdot\Omega,$$
where $\Omega=-\frac{x}{r^2}\wedge \nabla_x$ and $\wedge$ is the usual vector cross product on $\R^3$,
as well as the fact that
$$|\partial_r w(\cd)|\le |\nabla_x w(\cd)|,\quad \Bigl|\frac{1}{r}\Omega w(\cd)\Bigr|\lesssim |\nabla w(\cd)|.$$

\noindent{\em Proof of Lemma \ref{lemma.KS.estimates}:}
The proof will be broken into the cases that $|x|\le t/2$ and $|x|\ge t/2$.

\noindent{\em Case 1 ($|x|\le t/2$):}  Our main technique here will be the use of substitutions
involving $\Omega_{0i}=t\partial_i +x_i\partial_t$, as well as the fact that the commutator of any
$\Gamma$ with $\partial$ only involves $\partial$.

In order to obtain \eqref{KS.delta}, we begin with
\begin{align*}
\sum_{1\le i\le 3} \Omega_{0i}\partial_i w &= \sum_{1\le i\le 3} (t\partial_i + x_i\partial_t)\partial_i w\\
&= t \Delta w + \sum_{1\le i\le 3} x_i \partial_i \partial_t w\\
&= t\Delta w + \sum_{1\le i\le 3} \frac{x_i}{t} \Omega_{0i}\partial_t w - \frac{r^2}{t}\partial_t^2 w\\
&=\Bigl(t-\frac{r^2}{t}\Bigr)\Delta w + \sum_{1\le i\le 3} \frac{x_i}{t}\Omega_{0i}\partial_t w
-\frac{r^2}{t}\Box w.
\end{align*}
If we solve for the term involving $\Delta$ and use the fact that $t-\frac{r^2}{t}\ge t-r$ in this region,
we see that
$$(t-r)|\Delta w|\le \sum_{1\le i\le 3} |\Omega_{0i}\partial_i w|
+\sum_{1\le i\le 3} |\Omega_{0i}\partial_t w| + t|\Box w|,$$
from which \eqref{KS.delta} follows.

We use similar arguments to obtain \eqref{KS.dtdt}.  Indeed, by noting that
\begin{align*}
\sum_{1\le i\le 3} \Omega_{0i}\partial_i w &= t\Delta w + \sum_{1\le i\le 3} \frac{x_i}{t}\Omega_{0i}
\partial_t w - \frac{r^2}{t}\partial_t^2 w\\
&=\Bigl(t-\frac{r^2}{t}\Bigr)\partial_t^2 w + \sum_{1\le i\le 3} \frac{x_i}{t}\Omega_{0i} \partial_t w
-t\Box w,
\end{align*}
we see that \eqref{KS.dtdt} follows easily.

To prove \eqref{KS.dtdx}, we begin with 
\begin{align*}
(t-r)\partial_i\partial_t &= (t-r)\Bigl(\frac{1}{t}\Omega_{0i}-\frac{x_i}{t}\partial_t\Bigr)\partial_t\\
&=\Omega_{0i}\partial_t - x_i \partial_t^2 -\frac{r}{t}\Omega_{0i}\partial_t +\frac{rx_i}{t}\partial_t^2.
\end{align*}
Since we are in the region $r\le t/2$, which implies that $t\le \frac{1}{2}|t-r|$, this yields the bound
$$|(t-r)\partial_i \partial_t w|\lesssim |\Omega_{0i}\partial_t w|+|t-r||\partial_t^2 w|.$$
By applying \eqref{KS.dtdt}, we see that \eqref{KS.dtdx} follows.

\noindent{\em Case 2 ($|x|\ge t/2$):}
Here, we introduce the vector field
\begin{equation}
\label{Omega.r}
\Omega_r = \frac{x^i}{|x|}\Omega_{0i} = t\partial_r + r\partial_t,
\end{equation}
and begin by noticing that
\begin{align}
\partial_t \Omega_r w - \partial_r w &= t\partial_t \partial_r w + r\partial_t^2 w,\label{Omega.r.1}\\
\partial_r \Omega_r w - \partial_t w &= t\partial_r^2 w + r\partial_r \partial_t w.\label{Omega.r.2}
\end{align}
These immediately yield
\begin{equation}
\label{Omega.r.3}
r(\partial_t \Omega_r w - \partial_r w) - t(\partial_r \Omega_r w - \partial_t w) = r^2\partial_t^2 w 
-t^2 \partial_r^2 w.
\end{equation}

In order to prove \eqref{KS.delta}, we observe that \eqref{Omega.r.3} is equivalent to
$$r(\partial_t \Omega_r w-\partial_r w)-t(\partial_r\Omega_r w-\partial_t w) = r^2 \Box w
+t^2(\Delta w-\partial_r^2 w)-(t^2-r^2)\Delta w.$$
Rearranging terms, this is
$$(t-r)\Delta w = \frac{1}{t+r}\Bigl[r^2 \Box w+t^2(\Delta w-\partial_r^2 w) - r(\partial_t \Omega_r w
-\partial_r w)+t(\partial_r\Omega_r w-\partial_t w)\Bigr],$$
which yields the bound
$$|(t-r)\Delta w|\le r|\Box w|+t|\Delta w-\partial_r^2 w|+|\partial_t\Omega_r w|+|\partial_r w|
+|\partial_r \Omega_r w|+ |\partial_t w|.$$
Estimate \eqref{KS.delta} follows if we use \eqref{Delta.drdr} and that we are in the region $t\le 2r$.
Additionally, we require that
\begin{equation}
\label{Omega.r.d}
|\partial_t\Omega_r w| + |\partial_r \Omega_r w|\le \sum_{|\alpha|\le 1}|\partial \Gamma^\alpha w|,
\end{equation}
which follows directly from the definition \eqref{Omega.r}.

To prove \eqref{KS.dtdt}, we begin similarly by adding and subtracting $t^2\Box w$ to \eqref{Omega.r.3}
to get
$$r(\partial_t \Omega_r w-\partial_r w)-t(\partial_r \Omega_r w-\partial_t w)
=(r^2-t^2)\partial_t^2 w + t^2\Box w + t^2(\Delta w-\partial_r^2w).$$
This, in turn, yields
$$(t-r)\partial_t^2 w = \frac{1}{t+r}\Bigl[t(\partial_r \Omega_r w - \partial_t w)-r(\partial_t
\Omega_r w - \partial_r w)+t^2\Box w + t^2(\Delta w - \partial_r^2 w)\Bigr].$$
Using \eqref{Delta.drdr} and \eqref{Omega.r.d}, we see that \eqref{KS.dtdt} follows.

To prove \eqref{KS.dtdx}, we split the gradient into its radial and angular components,
$$\nabla_x = \frac{x}{r}\partial_r - \frac{x}{r^2}\wedge \Omega.$$
Beginning with the radial component, we see that
\begin{align*}
(t-r)\partial_t \partial_r &=\frac{t-r}{r}\Omega_r \partial_r - t\frac{t-r}{r}\partial_r^2\\
&=\frac{t-r}{r}\Omega_r\partial_r - t\frac{t-r}{r}\Bigl(\Delta - \frac{2}{r}\partial_r - \frac{1}{r^2}
\Omega\cdot\Omega\Bigr).
\end{align*}
Using that $t\le 2r$, we see that
\begin{align*}
|(t-r)\partial_t\partial_r w|&\le |\Omega_r \partial_r w|+|(t-r)\Delta w|+|\partial_r w|
+|\nabla \Omega w|\\
&\le |\partial_r\Omega_r w|+|\partial_t w|+|(t-r)\Delta w|+|\nabla w|+|\nabla\Omega w|.
\end{align*}
By \eqref{Omega.r.d} and \eqref{KS.delta}, we see that \eqref{KS.dtdx} holds when $\nabla_x$ is replaced
by $\partial_r$.

To prove \eqref{KS.dtdx} for the angular components of the gradient when $t\le 2r$, the bound follows trivially
as 
$$|(t-r)\frac{x}{r^2}\wedge \Omega \partial_t w|\lesssim |\partial_t \Omega w|,$$
which completes the proof.\qed

By arguing as in the proof of G\aa rding's inequality (see \cite{KSid}, Lemma 3.1), we can obtain
\begin{corr}
\label{corr.KS.L2}
Let $w\in C^\infty_0(\R_+\times\R^3)$.  Then,
\begin{equation}
\label{KS.L2}
\|\la t-r\ra \partial_{t,x}^2 \Gamma^\alpha w(t,\cd)\|_2\lesssim \sum_{|\beta|\le |\alpha|+1}
\|\Gamma^\beta w'(t,\cd)\|_2 + \sum_{|\beta|\le |\alpha|} \|\la t+r\ra \Box \Gamma^\beta u(t,\cd)\|_2.
\end{equation}
\end{corr}

This estimate will be used to control the last term in \eqref{hidano.1}, thus providing our main
wave equation decay estimate.

%%%%%%%%%%%%%%%%%%%%%%%%%%%%%%%%%%%%%%%%%%%%%%%%%%%%%%%%%%%%%%%%%%%%%%%%%%%%%%%%%%%%%%%%%%%%%%%%%%%%%
%%%%%%%%%%%%%%%%%%%%%%%%%%%%%%%%%%%%%%%%%%%%%%%%%%%%%%%%%%%%%%%%%%%%%%%%%%%%%%%%%%%%%%%%%%%%%%%%%%%%%
\newsection{Linear decay estimates in $\R_+\times\R^3$}

In this section, we gather the linear decay estimates that we shall require for the boundaryless
cases.  There are two main estimates.  The first is for the linear wave equation and will
be used to handle eigenfunctions with zero eigenvalue.  The latter is for linear Klein-Gordon
equations and will be used to estimate the eigenfunctions of nonzero modes.

Our decay estimate for the wave equation is a rather simple consequence of \eqref{sobolev.2}, \eqref{hidano.1},
and \eqref{KS.L2}.
\begin{proposition}
\label{prop.wave.decay}
Suppose that $w\in C^\infty(\R\times\R^3)$ satisfies $w(t,x)=0$, $t\le 2B$ where $B$ is a fixed
positive constant.  Suppose also that $\Box w(t,x)=0$ for $|x|>t-B$.  Then, 
\begin{equation}
\label{wave.decay}
(1+t)\sup_{x} |\partial_{t,x} w(t,x)|\lesssim \sum_{|\alpha|\le 2} \sum_k 
\sup_{\tau\in [2^{k-1},2^{k+1}]\cap [2B,t]} 2^k \|\Gamma^\alpha \Box w(\tau,\cd)\|_2.
\end{equation}
\end{proposition}
In this section, as in the previous, we are using $\Box=\partial_t^2-\Delta_{\R^3}$ as we are discussing 
the boundaryless wave equation.

\noindent{\em Proof of Proposition \ref{prop.wave.decay}:}
Here, again, we examine the regions $|x|\le t/2$ and $|x|>t/2$ separately.  In the former case, 
we may apply \eqref{hidano.1} and \eqref{KS.L2} to see that
$$(1+t)\sup_{\{x\,:\,|x|\le t/2\}} |\partial w(t,x)|\lesssim \sum_{|\alpha|\le 2} \|\Gamma^\alpha w'(t,\cd)\|_2
+ \sum_{|\alpha|\le 1} \|\la t+r\ra \Gamma^\alpha \Box w(t,\cd)\|_2.$$
By the energy inequality, this is
$$\lesssim \sum_{|\alpha|\le 2} \int_{2B}^t \|\Gamma^\alpha \Box w(\tau,\cd)\|_2\:d\tau
+ t\sum_{|\alpha|\le 1} \|\Gamma^\alpha \Box w(t,\cd)\|_2.$$
If we dyadically decompose $[2B,t]$, it follows immediately that this is
$$\lesssim \sum_{|\alpha|\le 2} \sum_k \sup_{\tau\in [2^{k-1},2^{k+1}]\cap[2B,t]}
2^k \|\Gamma^\alpha \Box w(\tau,\cd)\|_2$$
as desired.

Over $|x|>t/2$, the estimate follows similarly.  Indeed, we apply \eqref{sobolev.2} and the energy inequality
to see that 
$$(1+t)\sup_{\{x\,:\,|x|>t/2\}} |\partial w(t,x)| \lesssim \sum_{|\alpha|\le 2} \int_{2B}^t \|\Gamma^\alpha
\Box w(\tau,\cd)\|_2\:d\tau,$$
and use a dyadic decomposition in the temporal variable, as above, to see \eqref{wave.decay}.
\qed

The next estimate will be used to bound the eigenfunctions with non-vanishing eigenvalue.  This non-vanishing
eigenvalue will serve as a mass term after we expand, and thus, we are left with finding decay estimates
for Klein-Gordon equations.  As in \cite{MSS}, we will rely on the following estimate which is essentially
from H\"ormander \cite{H} (Proposition 7.3.6).
\begin{proposition}
\label{prop.kg.decay}
Suppose that $w\in C^\infty(\R\times\R^3)$ satisfies $w(t,x)=0$, $t\le 2B$ where $B$ is
a fixed positive constant.  Suppose also that $(\Box+\mu^2)w(t,x)=0$ for $|x|>t-B$.  Then,
\begin{equation}
\label{kg.decay}
(1+t)^{3/2} \sup_x |w(t,x)|\lesssim \sum_{|\alpha|\le 5}\sum_k
\sup_{\tau\in [2^{k-1},2^{k+1}]\cap [2B,t]} 2^k \|\Gamma^\alpha (\Box+\mu^2)w(\tau,\cd)\|_2
\end{equation}
where the implicit constant is independent of $\mu\ge 1$.
\end{proposition}

To prove this, one needs only modify the proof in \cite{H} by using the proper version of 
Lemma 7.3.4 of \cite{H}.  This states that if $v''+\mu^2 v=h$ in $[a,b]\subset \R$, then
$$\sup_{a\le \rho\le b} |v(\rho)|\le |v(a)|+|v'(a)|+\frac{1}{\mu}\int_a^b |h(\rho)|\:d\rho.$$

%%%%%%%%%%%%%%%%%%%%%%%%%%%%%%%%%%%%%%%%%%%%%%%%%%%%%%%%%%%%%%%%%%%%%%%%%%%%%%%%%%%%%%%%%%%%%%%%%%%%%%%
%%%%%%%%%%%%%%%%%%%%%%%%%%%%%%%%%%%%%%%%%%%%%%%%%%%%%%%%%%%%%%%%%%%%%%%%%%%%%%%%%%%%%%%%%%%%%%%%%%%%%%
\newsection{Linear decay estimates for waveguides}

In this section, we show that the estimates of the previous section can be adapted to yield
estimates in waveguides.  In the process, we lose some additional regularity which depends on the dimension
$d$ of $\Omega\subset \R^d$.  

Before we proceed, we review some basics from spectral theory and elliptic regularity theory.  
We refer the interested reader to, e.g., the
texts of Taylor \cite{Taylor} and Gilbarg and Trudinger \cite{gilbarg} for more thorough treatments.

Here, we allow $\Delta_\Omega$ to denote either the Dirichlet Laplacian with boundary conditions
\begin{equation}
\label{Dirichlet.bc}
h|_\bdy=0,
\end{equation}
or the Neumann Laplacian where the boundary conditions are
\begin{equation}\label{Neumann.bc}
\partial_\nu h|_\bdy=0,
\end{equation}
where $\partial_\nu$ denotes the normal derivative.  Since $\Omega$ is compact with smooth boundary,
it is known that the spectrum of $-\Delta_\Omega$ is discrete and nonnegative.  Letting
$\lambda_1^2\le \lambda_2^2 \le \lambda_3^2\le \dots$ denote the eigenvalues (counted with 
multiplicity), we know that
$$0<\lambda_1^2\le \lambda_2^2\le \dots\quad\text{for the Dirichlet Laplacian,}$$
and
$$0=\lambda_1^2<\lambda_2^2\le \lambda_3^2\le \dots\quad\text{for the Neumann Laplacian}.$$
We shall let $E_j : L^2(\Omega)\to L^2(\Omega)$ denote the projection onto the $j$th eigenspace.
Thus, for $h\in L^2(\Omega)$, we have that $E_jh$ is smooth and satisfies
$$-\Delta_\Omega E_j h(x)=\lambda_j^2 E_j h(x).$$
Moreover, we have a Plancherel's theorem
$$\|h\|^2_{L^2(\Omega)}=\sum_{j=1}^\infty \|E_jh\|_2^2.$$
Since by the Weyl formula, $\lambda_j\approx j^{1/d}$, $j=2,3,\dots$, we see that for $h\in C^\infty(\bar{\Omega})$,
\begin{equation}\label{weyl}
\begin{split}
(1+j)^{2/d} \|E_j h\|_{L^2(\Omega)} &\le C\|(I-\Delta_\Omega)E_j h\|_{L^2(\Omega)}\\
&=C \|E_j(I-\Delta_\Omega)h\|_{L^2(\Omega)},\quad j=1,2,3,\dots
\end{split}\end{equation}
provided that either \eqref{Dirichlet.bc} or \eqref{Neumann.bc} holds.

As for the elliptic regularity theory, we start with the basic estimate
\begin{lemma}
\label{lemma.ell.reg}
Suppose that $h\in C^\infty(\bar{\Omega})$ and that either \eqref{Dirichlet.bc} or \eqref{Neumann.bc} holds.
Then, for $N=2,3,\dots$,
\begin{equation}
\label{ell.reg}
\sum_{|\alpha|\le N} \|\partial^\alpha_y h\|_{L^2(\Omega)} \lesssim \sum_{|\alpha|\le N-2}
\|\partial^\alpha \Delta_\Omega h\|_{L^2(\Omega)} + \|h\|_{L^2(\Omega)}.
\end{equation}
\end{lemma}

By writing $-\Delta_\Omega=\Box-\partial_t^2+\Delta_{\R^3}$ and noticing that 
$\partial_t$ and $\nabla_x$ preserve the boundary conditions, it follows from \eqref{ell.reg}
and an induction argument that
\begin{corr}
\label{corr.ell.reg}
Suppose that $u(t,x,y)\in C^\infty(\R_+\times \R^3\times\Omega)$ and that either
$u(t,x,y)|_{y\in\bdy}=0$ or $\partial_\nu u(t,x,y)|_{y\in\bdy}=0$.  Then,
for $N=2,4,6,\dots$,
\begin{multline}
\label{ell.reg.2}
\sum_{|\alpha|\le N} \|\partial^\alpha_y u(t,x,\cd)\|_{L^2(\Omega)}\lesssim
\sum_{|\alpha|\le N} \|\partial_{t,x}^\alpha u(t,x,\cd)\|_{L^2(\Omega)}
\\+\sum_{|\alpha|\le N-2} \|\partial^\alpha_{t,x,y} \Box u(t,x,\cd)\|_{L^2(\Omega)},
\end{multline}
and for $N=2,3,4,\dots$,
\begin{multline}
\label{ell.reg.3}
\sum_{|\alpha|\le N} \|\partial^\alpha_y u(t,x,\cd)\|_{L^2(\Omega)}
\lesssim \sum_{\substack{|\alpha|+|\beta|\le N\\|\beta|\le 1}} \|\partial^\alpha_{t,x}\partial_y^\beta
u(t,x,\cd)\|_{L^2(\Omega)} \\+ \sum_{|\alpha|\le N-2} \|\partial^\alpha_{t,x,y} \Box u(t,x,\cd)\|_{L^2(\Omega)}.
\end{multline}
\end{corr}

The following lemma shows how our equation behaves under the eigenfunction expansion.
\begin{proposition}
\label{prop.efunc.exp}
Assume that $u(t,x,y)\in C^\infty(\R_+\times\R^3\times\Omega)$ satisfies either $u(t,x,y)|_{y\in\bdy}=0$
or $\partial_\nu u(t,x,y)|_{y\in\bdy}=0$.  Then, it follows that
$$E_j \Box u(t,x,y)=(\partial_t^2-\Delta_{\R^3}+\lambda_j^2)E_ju(t,x,y)$$
and
$$\sum_j \|(\partial_t^2-\Delta_{\R^3}+\lambda_j^2)E_ju(t,x,\cd)\|^2_{L^2(\Omega)}
=\|\Box u(t,x,\cd)\|^2_{L^2(\Omega)}.$$
\end{proposition}

The last of our basic elliptic regularity estimates follows from Sobolev's lemma and
\eqref{ell.reg.2}.
\begin{proposition}
\label{prop.ell.reg}
Assume that $u(t,x,y)\in C^\infty(\R_+\times \R^3\times\Omega)$ and that either $u(t,x,y)|_{y\in\bdy}=0$
or $\partial_\nu u(t,x,y)|_{y\in \bdy}=0$.  Then,
\begin{equation}
\label{ell.reg.sobolev}
|\partial_y^\beta u(t,x,y)|\lesssim \sum_{|\alpha|\le |\beta|+(d+4)/2} \|\partial^\alpha_{t,x} 
u(t,x,\cd)\|_{L^2(\Omega)}
+ \sum_{|\alpha|\le |\beta|+d/2} \|\partial^\alpha_{t,x,y}\Box u(t,x,\cd)\|_{L^2(\Omega)}.
\end{equation}
\end{proposition}

We are now ready to present our main decay estimates.  We let
\begin{equation}
\label{tGamma}
\{\tGamma\}=\{\Gamma\}\cap \{\partial_y\}=\{\partial_t,\partial_x, \Omega_{jk},
\partial_y\,:\, 0\le j<k\le 3\}
\end{equation}
denote the full set of ``admissible'' vector fields.
\begin{proposition}
\label{proposition.main.decay}
Fix $B$ and suppose that $u\in C^\infty(\R_+\times\R^3\times\bar{\Omega})$ satisfies $u(t,x,y)=0$
for $t\le 2B$, $\Box u(t,x,y)=0$ for $|x|>t-B$, and $\partial_\nu
u(t,x,y)|_{y\in\bdy} = 0$.  Then,
\begin{multline}
\label{main.decay}
(1+t)|\tGamma^\beta \partial_{t,x} u(t,x,y)|\lesssim \sum_{|\alpha|\le |\beta|+(5d+16)/2}
\sum_k \sup_{\tau\in [2^{k-1},2^{k+1}]\cap [2B,t]} 2^k \|\tGamma^\alpha \Box u(\tau,\cd)\|_2
\\+(1+t) \sum_{|\alpha|\le |\beta|+(5d+6)/2} \|\tGamma^\alpha \Box u(t,\cd)\|_2.
\end{multline}
Here, as above, $d$ is the dimension of $\Omega$.
\end{proposition}

\noindent{\em Proof of Proposition \ref{proposition.main.decay}:}
The proof of this proposition follows that of Proposition 3.5 of \cite{MSS} very closely.  Indeed, the only
major modification is the use of \eqref{wave.decay} to estimate those terms in our expansion 
with vanishing eigenvalues.

By noting that $\Gamma$ preserves the Neumann boundary condition and commutes with $\Box$, it suffices
to take $\tGamma^\beta = \partial_y^\beta$.  Using \eqref{ell.reg.sobolev} and Plancherel's theorem, we have
that
\begin{equation}\label{step1} \begin{split}
|\partial^\beta_y \partial_{t,x} u(t,x,y)|^2 &\lesssim \sum_{|\alpha|\le |\beta|+(d+4)/2}
\|\partial^\alpha_{t,x} \partial_{t,x}u(t,x,\cd)\|^2_{L^2(\Omega)}
\\&\qquad\qquad\qquad\qquad\qquad\qquad
+\sum_{|\alpha|\le |\beta|+(d+2)/2} \|\partial^\alpha_{t,x,y}\Box u(t,x,\cd)\|^2_{L^2(\Omega)}\\
&\lesssim \sum_{|\alpha|\le |\beta|+(d+4)/2} \sum_{j=1}^\infty \|E_j \partial^\alpha_{t,x} \partial_{t,x}
u(t,x,\cd)\|^2_{L^2(\Omega)}\\
&\qquad\qquad\qquad\qquad\qquad\qquad
+\sum_{|\alpha|\le |\beta|+(d+6)/2} \|\partial^\alpha_{t,x,y}\Box u(t,\cd)\|^2_{L^2(\R^3\times\Omega)}.
\end{split}\end{equation}
In the last step, we have also applied Sobolev's lemma.

In order to gain further control over the first term in the right, we apply \eqref{weyl} to see that
\begin{align*}
(1+j)^{2/d}& \sum_{|\alpha|\le |\beta|+(d+4)/2} \|E_j\partial^\alpha_{t,x} \partial_{t,x} u(t,x,\cd)\|_{L^2(\Omega)}\\
&\lesssim \sum_{|\alpha|\le |\beta|+(d+4)/2} \|E_j (I-\Delta_\Omega)\partial^\alpha_{t,x}\partial_{t,x}u(t,x,
\cd)\|_{L^2(\Omega)}\\
&\lesssim \sum_{|\alpha|\le |\beta|+(d+8)/2} \|E_j \partial^\alpha_{t,x}\partial_{t,x} u(t,x,\cd)\|_{L^2(\Omega)}\\
&\qquad\qquad\qquad\qquad\qquad\qquad\qquad
+\sum_{|\alpha|\le |\beta|+(d+6)/2} \|E_j\partial^\alpha_{t,x} \Box u(t,x,\cd)\|_{L^2(\Omega)}.
\end{align*}
In the last step, we have simply written $-\Delta_\Omega=\Box -\partial_t^2 +\Delta$.  If we apply
this estimate recursively, it follows that
\begin{multline*}
\sum_{|\alpha|\le |\beta|+(d+4)/2} \|E_j \partial_{t,x}^\alpha \partial_{t,x}u(t,x,\cd)\|_{L^2(\Omega)}
\\\lesssim \frac{1}{(1+j)^2} \sum_{|\alpha|\le |\beta|+(5d+4)/2} \|E_j \partial_{t,x}^\alpha \partial_{t,x}
u(t,x,\cd)\|_{L^2(\Omega)}
\\+\sum_{|\alpha|\le |\beta|+(5d+2)/2} \|E_j \partial^\alpha_{t,x}\Box u(t,x,\cd)\|_{L^2(\Omega)}.
\end{multline*}
By Plancherel's theorem and Sobolev's lemma on $\R^3$, this implies that
\begin{multline}\label{step2}
\sum_{|\alpha|\le |\beta|+(d+4)/2} \sum_{j=1}^\infty \|E_j \partial_{t,x}^\alpha \partial_{t,x} u(t,x,\cd)\|_{
L^2(\Omega)}^2
\lesssim \sum_{|\alpha|\le |\beta|+(5d+6)/2} \|\partial^\alpha_{t,x} \Box u(t,\cd)\|^2_{L^2(\R^3\times\Omega)}
\\+\sum_{j=1}^\infty \frac{1}{(1+j)^4}\sum_{|\alpha|\le |\beta|+(5d+4)/2} \|E_j \partial^\alpha_{t,x}
\partial_{t,x} u(t,x,\cd)\|^2_{L^2(\Omega)}.
\end{multline}

Combining \eqref{step1} and \eqref{step2}, it follows that
\begin{multline}
\label{step3}
(1+t)^2 |\partial^\beta_y \partial_{t,x}u(t,x,y)|^2 \lesssim (1+t)^2\sum_{|\alpha|\le |\beta|+(5d+6)/2}
\|\partial^\alpha_{t,x} \Box u(t,\cd)\|^2_{L^2(\R^3\times\Omega)}
\\
+(1+t)^2 \sum_{|\alpha|\le |\beta|+(5d+4)/2} \|\partial_{t,x} E_1 \partial^\alpha_{t,x} u(t,x,\cd)\|^2_{L^2(\Omega)}
\\+(1+t)^2 \sum_{j=2}^\infty \frac{1}{(1+j)^4} \sum_{|\alpha|\le |\beta|+(5d+6)/2} \|
E_j \partial^\alpha_{t,x}
 u(t,x,\cd)\|^2_{L^2(\Omega)}.
\end{multline}

Here, we apply \eqref{wave.decay} to the second term on the right side, which yields
\begin{align*}
(1+t)&\sum_{|\alpha|\le |\beta|+(5d+4)/2} \|\partial_{t,x}E_1\partial_{t,x}^\alpha u(t,x,\cd)\|_{L^2(\Omega)}
\\&\lesssim \sum_{|\alpha|\le |\beta|+(5d+8)/2} \sum_k \sup_{\tau\in [2^{k-1},2^{k+1}]\cap [2B,t]}
2^k \|\Gamma^\alpha (\partial_t^2-\Delta)E_1u(\tau,\cd)\|_{L^2(\R^3\times\Omega)}\\
&\lesssim \sum_{|\alpha|\le |\beta|+(5d+8)/2} \sum_k \sup_{\tau\in [2^{k-1},2^{k+1}]\cap [2B,t]}
2^k \|\Gamma^\alpha \Box u(\tau,\cd)\|_{L^2(\R^3\times\Omega)}
\end{align*}
The last inequality follows from Proposition \ref{prop.efunc.exp}.

We can similarly apply \eqref{kg.decay} to each of the summands in the last term of \eqref{step3},
where $\lambda_j\ge c>0$ for each $j=2,3,\dots$.  This yields
\begin{align*}
(1+t) &\sum_{|\alpha|\le |\beta|+(5d+6)/2} \|E_j \partial^\alpha_{t,x}u(t,x,\cd)\|_{L^2(\Omega)}
\\&\lesssim \sum_{|\alpha|\le |\beta|+(5d+16)/2} \sum_k \sup_{\tau\in [2^{k-1},2^{k+1}]\cap [2B,t]}
2^k \|\Gamma^\alpha (\partial_t^2-\Delta+\lambda_j^2)E_j u(\tau,\cd)\|_{L^2(\R^3\times\Omega)}\\
&\lesssim \sum_{|\alpha|\le |\beta|+(5d+16)/2} \sum_k \sup_{\tau\in [2^{k-1},2^{k+1}]\cap [2B,t]}
2^k \|\Gamma^\alpha \Box u(\tau,\cd)\|_{L^2(\R^3\times\Omega)}.
\end{align*}
Here, again, the last inequality is a result of Proposition \ref{prop.efunc.exp}.

Combining these last two bounds with \eqref{step3} completes the proof.\qed

For the case of Dirichlet boundary conditions, a similar bound can be achieved.
\begin{proposition}
\label{prop.dirichlet.decay}
Fix $B$ and suppose that $u\in C^\infty(\R_+\times \R^3\times\bar{\Omega})$ satisfies $u(t,x,y)=0$
for $t\le 2B$, $\Box u(t,x,y)=0$ for $|x|>t-B$, and $u(t,x,y)|_{y\in\bdy}=0$.  Then,
\begin{multline}
\label{dirichlet.decay}
(1+t) |\tGamma^\beta u(t,x,y)|\lesssim \sum_{|\alpha|\le |\beta|+(5d+14)/2}
\sum_k \sup_{\tau\in [2^{k-1},2^{k+1}]\cap [2B,t]} 2^k \|\tGamma^\alpha \Box u(\tau,\cd)\|_2
\\+ (1+t) \sum_{|\alpha|\le |\beta|+(5d+5)/2} \|\tGamma^\alpha \Box u(t,\cd)\|_2.
\end{multline}
\end{proposition}
This is essentially Proposition 3.5 from \cite{MSS}.  In fact, as there is no zero mode, an estimate
with stronger decay is available.  In our proof of long time existence, however, we would not be
able to handle the added power of
$(1+t)^{1/2}$ that would be necessary in the last term of \eqref{dirichlet.decay}.  The modification
required to obtain this weakened version is straightforward, and as the proof closely resembles that of
the previous proposition, it will be omitted.

Using this, we can obtain the necessary decay for $\partial_y u$ when $\Omega$ is one-dimensional.  Indeed,
we simply notice that when $\Omega$ is one-dimensional, $v=\partial_y u = \partial_\nu u$ satisfies 
the Dirichlet wave equation $\Box v = \partial_y \Box u$.  Thus, we may apply \eqref{dirichlet.decay}
to obtain the following.
\begin{corr}
\label{corr.dy.decay}
Fix $B$ and suppose that $u\in C^\infty(\R_+\times\R^3\times [a,b])$ for two fixed constants $a<b$.  
Suppose further that $u$ satisfies $u(t,x,y)=0$ for $t\le 2B$, $\Box u(t,x,y)=0$ for $|x|>t-B$, and
$\partial_y u(t,x,a)=\partial_y u(t,x,b)=0$.  Then,
\begin{multline}
\label{dy.decay}
(1+t)|\tGamma^\beta \partial_y u(t,x,y)|\lesssim \sum_{|\alpha|\le |\beta|+10}
\sum_k \sup_{\tau\in [2^{k-1},2^{k+1}]\cap [2B,t]} 2^k \|\tGamma^\alpha \Box u(\tau,\cd)\|_2
\\+(1+t)\sum_{|\alpha|\le |\beta|+6} \|\tGamma^\alpha \Box u(t,\cd)\|_2.
\end{multline}
\end{corr}

%%%%%%%%%%%%%%%%%%%%%%%%%%%%%%%%%%%%%%%%%%%%%%%%%%%%%%%%%%%%%%%%%%%%%%%%%%%%%%%%%%%%%%%%%%%%%%%%%%%%%
%%%%%%%%%%%%%%%%%%%%%%%%%%%%%%%%%%%%%%%%%%%%%%%%%%%%%%%%%%%%%%%%%%%%%%%%%%%%%%%%%%%%%%%%%%%%%%%%%%%%%
\newsection{Energy estimates}

In this section, we briefly describe the energy estimates for the variable coefficient wave equation
which we shall require.  For $\gamma^{jk}\in C^\infty$ satisfying
\begin{equation}\label{perturbation.smallness}
\sum_{j,k=0}^{3+d} |\gamma^{jk}| \le 1/2,\quad
\sum_{i,j,k=0}^{3+d} \|\partial_i \gamma^{jk}\|_{L^1_tL^\infty_{x,y}([0,T]\times \R^3\times\Omega)}
\le C_0,
\end{equation}
we study
$$\begin{cases}
  \Box u + \sum_{j,k=0}^{3+d} \gamma^{jk}(t,x,y)\partial_j\partial_k u = F,\quad 2B\le t\le T,\\
  \partial_\nu w(t,x,y)=0,\quad y\in \bdy,\\
  w(t,x,y)=0,\quad t\le 2B.
  \end{cases}
$$
In order to establish energy estimates with Neumann boundary conditions, we are required to assume
a natural nonlinear compatibility condition on the $\gamma^{jk}$.  This says that
\begin{equation}\label{Neumann.assumption}\begin{split}
\sum_{j,k=0}^{3+d} \gamma^{jk}(t,x,y)&\xi_j\theta_k = 0,\quad \text{if }y\in\bdy,\\
\theta=(0,\dots,0,\nu_1(y),\dots,&\nu_d(y)),\quad \text{ and } \xi\cdot \theta=0.
\end{split}\end{equation}
Under this assumption, if $u$ vanishes for large $|x|$, 
the standard proof of energy inequalities (see, e.g., Proposition 2.1 in Chapter
1 of \cite{Sogge}) yields
\begin{equation}
\label{basic.energy}
\|\partial_{t,x,y} u(t,\cd)\|_2 \lesssim \int_0^t \|F(s,\cd)\|_2\:ds. 
\end{equation}

In the sequel, we shall require a version of \eqref{basic.energy} where $u$ is replaced by $\tGamma^\beta u$.
Since $\Gamma$ preserves the boundary condition, it follows from \eqref{ell.reg.3} that
$$\sum_{|\alpha|+|\beta|\le N}\|\Gamma^\alpha \partial_y^\beta \partial_{t,x,y}u(t,\cd)\|_2
\lesssim \sum_{|\alpha|\le N} \|\Gamma^\alpha \partial_{t,x,y} u(t,\cd)\|_2
+\sum_{|\alpha|\le N-1} \|\tGamma^\alpha \Box u(t,\cd)\|_2.$$
Thus, since
$$\Box \Gamma^\alpha u + \sum_{j,k=0}^{3+d} \gamma^{jk}\partial_j\partial_k \Gamma^\alpha u
=\Gamma^\alpha F + \sum_{j,k=0}^{3+d} [\gamma^{jk},\Gamma^\alpha] \partial_j\partial_k w
+\sum_{j,k=0}^{3+d} \gamma^{jk}[\partial_j\partial_k, \Gamma^\alpha]w,$$
it follows easily that for $N=0,1,2,\dots$
\begin{multline}
\label{main.energy}
\sum_{|\alpha|\le N} \|\tGamma^\alpha \partial_{t,x,y} u(t,\cd)\|_2
\lesssim \sum_{|\alpha|\le N} \int_0^t \|\Gamma^\alpha F(s,\cd)\|_2\:ds
\\+ \sum_{|\alpha|\le N} \sum_{j,k=0}^{3+d} \int_0^t \|[\gamma^{jk},\Gamma^\alpha]\partial_j \partial_k u(s,\cd)\|_2
\:ds
+\sum_{|\alpha|\le N} \sum_{j,k=0}^{3+d} \int_0^t \|\gamma^{jk} [\Gamma^\alpha, \partial_j\partial_k]u(s,\cd)\|_2\:ds
\\+\sum_{|\beta|\le N-1} \|\tGamma^\beta \Box u(t,\cd)\|_2.
\end{multline}

%%%%%%%%%%%%%%%%%%%%%%%%%%%%%%%%%%%%%%%%%%%%%%%%%%%%%%%%%%%%%%%%%%%%%%%%%%%%%%%%%%%%%%%%%%%%%%%%%%%%%
%%%%%%%%%%%%%%%%%%%%%%%%%%%%%%%%%%%%%%%%%%%%%%%%%%%%%%%%%%%%%%%%%%%%%%%%%%%%%%%%%%%%%%%%%%%%%%%%%%%%%
\newsection{Almost global existence}

In this section, we prove Theorem \ref{theorem1.1} and Theorem \ref{theorem1.2}.  We will prove 
the former explicitly and will mention in the course of this proof the modifications that would yield
the latter.  In order to simplify the exposition, we will take $d=1$ in Theorem \ref{theorem1.1} also.
The general case only requires additional occurrences of the vector fields (the number depending on $d$) 
in the definitions \eqref{Mk} and \eqref{Ak} below.  

In order to apply our decay estimates, we shift the time variable by $2B$ and examine initial conditions
at time $t=2B$.  Using local existence theory (see, e.g., \cite{KSS}) if $\varepsilon>0$ is sufficiently small and
$N$ is sufficiently large in \eqref{data.smallness}, it is known that there is a solution to \eqref{main.equation}
on an interval $[2B, 2B+1]$.  We will use this to reduce to the case of vanishing data (as was similarly
done in \cite{KSS2}).  Fixing $\eta\in C^\infty(\R)$ with $\eta(t)\equiv 1$ for $t\le 2B+\frac{1}{2}$ and
$\eta(t)\equiv 0$ for $t\ge 2B+1$, we set
$$u_0(t,x,y)=\eta(t)u(t,x,y)$$
and note that this solves
$$\Box u_0 = \eta Q(\partial_{t,x}u,\partial_{t,x}^2 u)+[\Box,\eta]u$$
with 
$$\partial_\nu u_0(t,x,y)|_{y\in\bdy}=0.$$
We also note that the local existence results imply that
\begin{equation}
\label{local.small}
\sup_t \sum_{|\alpha|\le N} \|\tGamma^\alpha u_0(t,\cd)\|_{L^2(\R^3\times\Omega)}\le C_0 \varepsilon.
\end{equation}
We conclude that $u$ solves \eqref{main.equation} for $2B<t<T$ if and only if $w=u-u_0$ solves
\begin{equation}
\label{vanishing.data.equation}
\begin{cases}
\Box w= (1-\eta) Q(\partial_{t,x}(u_0+w),\partial_{t,x}^2(u_0+w))-[\Box,\eta](u_0+w),\\
\partial_\nu w(t,x,y)=0,\quad y\in \bdy,\\
w(t,x,y)=0,\quad t\le 2B.
\end{cases}
\end{equation}

We solve \eqref{vanishing.data.equation} on $[0,\T)$ using an iteration.  We set $w_0\equiv 0$ and
define $w_k$ recursively to solve
\begin{equation}
\label{wk.equation}
\begin{cases}
\Box w_k = (1-\eta)Q(\partial_{t,x}(u_0+w_{k-1}),\partial_{t,x}^2 (u_0+w_k))-[\Box,\eta](u_0+w_k),\\
\partial_\nu w_k(t,x,y)=0,\quad y\in\bdy,\\
w_k(t,x,y)=0,\quad t\le 2B,
\end{cases}
\end{equation}
for $k=1,2,3,\dots$.

We define
\begin{equation}
\label{Mk}
M_k(T)=\sup_{2B\le t\le T} \Bigl(\sum_{|\alpha|\le 30} \|\tGamma^\alpha \partial_{t,x,y} w_k(t,\cd)\|_2
+ \sum_{|\beta|\le 16} (1+t)\sup_{x,y} |\tGamma^\beta \partial_{t,x} w_k(t,x,y)|\Bigr).
\end{equation}
Here, in the proof of Theorem \ref{theorem1.2}, we need to, of course, allow $\partial_y w_k$ in the last
term.  We will label the two terms in \eqref{Mk} $I_k(T)$ and $II_k(T)$ respectively.

Our first goal is to inductively show that for $\varepsilon$ sufficiently small, we have
\begin{equation}
\label{Mk.goal}
M_k(\T)\le 4C_1\varepsilon.
\end{equation}
Here, $C_1$ is a uniform constant that is, say, 10 times greater than $(1+C_0)$ times the
implicit constants appearing in \eqref{main.decay} and \eqref{main.energy}.  When $k=1$,
\eqref{Mk.goal} follows from \eqref{data.smallness}, the well-known local estimates, \eqref{main.energy},
\eqref{main.decay}, and Gronwall's inequality.

Assuming that \eqref{Mk.goal} holds for $k$ replaced by $k-1$, we apply \eqref{main.energy} to see that
\begin{multline}
\label{Ik}
I_k(\T)\le C\int_0^\T \sum_{|\beta|\le 15} \|\tGamma^\beta \partial_{t,x} (u_0+w_{k-1})(t,\cd)\|_\infty
\sum_{|\alpha|\le 30} \|\tGamma^\alpha \partial_{t,x} (u_0+w_k)(t,\cd)\|_2\:dt
\\+C\int_0^\T \sum_{|\beta|\le 15} \|\tGamma^\beta \partial_{t,x}(u_0+w_k)(t,\cd)\|_\infty
\sum_{|\alpha|\le 30} \|\tGamma^\alpha \partial_{t,x}(u_0+w_{k-1})(t,\cd)\|_2\:dt\\
\\+C\int_0^\T \sum_{|\beta|\le 15} \|\tGamma^\beta \partial_{t,x}(u_0+w_{k-1})(t,\cd)\|_\infty
\sum_{|\alpha|\le 30} \|\tGamma^\alpha \partial_{t,x}(u_0+w_{k-1})(t,\cd)\|_2\:dt\\
+C\sup_{0\le t\le \T} \sum_{|\beta|\le 15} \|\tGamma^\beta \partial_{t,x} (u_0+w_{k-1})(t,\cd)\|_\infty
\sum_{|\alpha|\le 30} \|\tGamma^\alpha \partial_{t,x}(u_0+w_k)(t,\cd)\|_2\\
+C\sup_{0\le t\le \T} \sum_{|\beta|\le 15} \|\tGamma^\beta \partial_{t,x} (u_0+w_k)(t,\cd)\|_\infty
\sum_{|\alpha|\le 30} \|\tGamma^\alpha \partial_{t,x}(u_0+w_{k-1})(t,\cd)\|_2
\\+C\sup_{0\le t\le\T} \sum_{|\beta|\le 15} \|\tGamma^\beta \partial_{t,x} (u_0+w_{k-1})(t,\cd)\|_\infty
\sum_{|\alpha|\le 30} \|\tGamma^\alpha \partial_{t,x}(u_0+w_{k-1})(t,\cd)\|_2
\\+C_1\varepsilon + \sum_{|\alpha|\le 30} \int_{2B+(1/2)}^{2B+1} \|\Gamma^\alpha \partial_t w_k(t,\cd)\|_2\:dt. 
\end{multline}
Here, we have set
$$\gamma^{jk}=-\sum_{l=0}^3 A^{jk}_l\partial_l w_{k-1}.$$
Notice that since $\gamma^{jk}$ vanishes if either $j=n+1$ or $k=n+1$ (recall that we are looking
at $d=1$), \eqref{Neumann.assumption} is trivially satisfied.  For the proof of Theorem \ref{theorem1.2},
we instead need to refer to \eqref{Neumann.compatibility} to see that \eqref{Neumann.assumption} holds.

By using \eqref{Mk}, \eqref{local.small}, the inductive hypothesis, and Gronwall's inequality, 
we see that
\begin{equation}\begin{split}
\label{Ik.2}
I_k(\T)&\le C_1\varepsilon + C\varepsilon^2 + C\varepsilon M_{k-1}(\T) + C\varepsilon M_k(\T)
%+C(M_{k-1}(\T))^2
\\&\qquad\qquad\qquad\qquad\qquad
+C [(M_{k-1}(\T))^2+M_{k-1}(\T)M_k(\T)]\Bigl(1+\int_{2B}^\T \frac{1}{1+t}\:dt\Bigr)\\
&\le C_1\varepsilon + C\varepsilon^2 + 16CC_1^2\varepsilon^2(1+\log(1+\T)) + C\varepsilon M_k(\T)
\\
&\qquad\qquad\qquad\qquad\qquad\qquad\qquad\qquad\qquad\qquad+ 4CC_1 \varepsilon (1+\log(1+\T)) M_k(\T).
\end{split}
\end{equation}

In order to bound $II_k(\T)$, we apply \eqref{main.decay} (and notice that the second term in the right
side of \eqref{main.decay} is controlled by the preceding term).  In the proof of Theorem \ref{theorem1.2},
we would need to also use \eqref{dy.decay}.  This yields
\begin{multline}
\label{IIk}
II_k(\T)\\
\le \sum_j \sup_{\tau\in [2^{j-1},2^{j+1}]\cap [2B,\T]} 2^j \sum_{|\beta|\le 13} \|\tGamma^\beta
\partial_{t,x} (u_0+w_{k-1})(\tau,\cd)\|_\infty \sum_{|\alpha|\le 26} \|\tGamma^\alpha \partial_{t,x}
(u_0+w_k)(\tau,\cd)\|_2\\
+\sum_j \sup_{\tau\in [2^{j-1},2^{j+1}]\cap [2B,\T]} 2^j \sum_{|\beta|\le 13} \|\tGamma^\beta
\partial_{t,x} (u_0+w_k)(\tau,\cd)\|_\infty \sum_{|\alpha|\le 26} \|\tGamma^\alpha \partial_{t,x}
(u_0+w_{k-1})(\tau,\cd)\|_2\\
+\sum_j \sup_{\tau\in [2^{j-1},2^{j+1}]\cap [2B,\T]} 2^j \sum_{|\beta|\le 13} \|\tGamma^\beta
\partial_{t,x}(u_0+w_{k-1})(\tau,\cd)\|_\infty \sum_{|\alpha|\le 26} \|\tGamma^\alpha \partial_{t,x}
(u_0+w_{k-1})(\tau,\cd)\|_2\\
+C_1\varepsilon + C \sup_{2B+(1/2)\le t\le 2B+1} \sum_{|\alpha|\le 26} \|\tGamma^\alpha \partial_t
w_k(t,\cd)\|_2.
\end{multline}
Here, we use \eqref{Mk}, \eqref{local.small}, and \eqref{Ik.2} to see that
\begin{multline*}
II_k(\T)\le 2C_1\varepsilon + C\varepsilon^2 + C\varepsilon M_{k-1}(\T) + C\varepsilon M_k(\T)
\\+ C\Bigl(1+\log(2+\T)\Bigr)[M_k(\T)M_{k-1}(\T)+(M_{k-1}(\T))^2].\end{multline*}
In addition to that occurring in \eqref{Ik.2}, 
the $\log(2+\T)$ term appears due to the simple observation that we are summing over $O(\log(2+\T))$
choices of $j$ with $2^{j-1}\le \T$.
Applying the inductive hypothesis, we see that this is
\begin{multline}
\label{IIk.2}
II_k(\T)\le 2C_1\varepsilon + C\varepsilon^2 + 16CC_1^2\varepsilon^2(1+\log(2+\T)) + C\varepsilon M_k(\T)
\\+ 4CC_1 \varepsilon (1+\log(2+\T))M_k(\T).
\end{multline}

If we combine \eqref{Ik.2} and \eqref{IIk.2}, we see that this yields \eqref{Mk.goal} provided that
$\varepsilon$ is sufficiently small and $\T$ is as in \eqref{lifespan} with $\kappa$ sufficiently small.

It remains to show that the sequence $w_k$ converges for $t\in [2B,\T)$.  By setting
\begin{multline}
\label{Ak}
A_k(T)=\sup_{2B\le t\le T} \Bigl(\sum_{|\alpha|\le 30} \|\tGamma^\alpha \partial_{t,x} (w_k-w_{k-1})(t,\cd)\|_2
\\+ \sum_{|\beta|\le 16} (1+t) \sup_{x,y} |\tGamma^\beta \partial_{t,x} (w_k-w_{k-1})(t,x,y)|\Bigr)
\end{multline}
and arguing as above, it can be shown that
$$A_k(\T)\le \frac{1}{2}A_{k-1}(\T).$$
This suffices to show convergence and thus finishes the proof.\qed

%%%%%%%%%%%%%%%%%%%%%%%%%%%%%%%%%%%%%%%%%%%%%%%%%%%%%%%%%%%%%%%%%%%%%%%%%%%%%%%%%%%%%%%%%%%%%%%%%%%%%
%%%%%%%%%%%%%%%%%%%%%%%%%%%%%%%%%%%%%%%%%%%%%%%%%%%%%%%%%%%%%%%%%%%%%%%%%%%%%%%%%%%%%%%%%%%%%%%%%%%%%

\end{document}